\documentclass[graybox]{svmult}
\usepackage{mathptmx}       
\usepackage{helvet}         
\usepackage{courier}        
\usepackage{type1cm}        
\usepackage{makeidx}         
\usepackage{graphicx}        
\usepackage{multicol}        
\usepackage[bottom]{footmisc}
\usepackage{url}

\usepackage{amsfonts}
\usepackage{amsmath}
\usepackage{amssymb}

\begin{document}

\title{Deviation Identity for Linear Differential Operators and Its Application to Obstacle Problems}
\titlerunning{Deviation Identity for LDO and in Obstacle Problems} 
\author{Kseniya Darovskaya} 
\institute{Kseniya Darovskaya \at RUDN University, 
117198, Russia, Moscow,  Miklukho-Maklaya str.~6 \at \email{k.darovsk@gmail.com} } 

\maketitle

\abstract{For a class of variational problems with linear differential operator, we obtain a~convenient form of the deviation identity, i.e., the value of the distance between approximate solutions and the exact ones.  
We illustrate the result with an explicit form of the deviation identity for a biharmonic obstacle problem.} 



\section{Deviation identity}
Let us consider a linear bounded operator $\Lambda$ acting from $V,$ that is, a Sobolev space $H^k(\Omega)$ with some boundary conditions, into~$Y=Y(\Omega).$ It satisfies the coercivity condition 
\begin{equation*}\label{dka_сoer}
	\|\Lambda v \|_Y \ge \kappa \| v\|_V \quad (\forall v \in \mathbb K). 
\end{equation*}

Assume that the domain $\Omega \subset \mathbb R^n$ is bounded with Lipschitz boundary~$\partial\Omega,$ $Y$ is a space of scalar, vector or symmetric matrix $L_2$--functions with standard scalar products, and $\mathbb K \subseteq V$ is a nonempty convex closed subset.

The tensor $\mathcal A$ with elements from $L_\infty (\Omega)$ provides a linear continuous  self--adjoint bounded mapping $Y \to Y$ and satisfies the inequality  
\begin{equation*}\label{dka_condA}
	\kappa_1 \|\sigma\|^2_Y \le \left( \mathcal A \sigma , \sigma \right)_Y \le \kappa_2 \|\sigma\|^2_Y \quad (\forall \sigma \in Y),
\end{equation*}
where the positive constants $\kappa_1$ and $\kappa_2$ do not depend on $\sigma.$
It defines the equivalent norm $\| \gamma \|_{\mathcal A}=(\mathcal A \gamma, \gamma)^{{1}/{2}}_Y.$

Let $\mathcal A^{-1}$ be a mapping dual to~$\mathcal A,$ with the same properties, and let it satisfy the inequality 
\begin{equation*}\label{dka_condA-1}
	\kappa_2^{-1} \|\sigma^*\|^2_{Y^*} \le 
	\left( \mathcal A^{-1} \sigma^* , \sigma^* \right)_{Y^*}
	\le \kappa_1^{-1} \|\sigma^*\|^2_{Y^*} \quad (\forall \sigma^* \in Y^*).
\end{equation*} 
It also defines the equivalent norm: $\| \gamma^{\ *} \|_{\mathcal A^{-1}}=(\mathcal A^{-1} \gamma^{\ *}, \gamma^{\ *})^{{1}/{2}}_{Y^*}.$
 
Everywhere else we will omit lower indices in standard $L_2$ scalar products and norms. (Since both $Y$ and $Y^*$ are Lebesgue spaces, the last remark refers to them as well.)

Let $f \in L_2(\Omega).$ We now consider the following pair of variational problems:
\newline the primal one,
\begin{equation} \label{dka_primal_problem}
\frac{1}{2} \| \Lambda v \|^2_{\mathcal A} -  \left( f,v \right) \xrightarrow {} \min_{\mathbb K}
\end{equation}
\newline and the dual one,
\begin{equation} \label{dka_dual_problem}
	-\frac{1}{2} \| y^* \|^2_{\mathcal A^{-1}} - \sup\limits_{v \in \mathbb K} \biggl\{ - \left( y^*, \Lambda v \right) + (f,v) \biggr\}
	\xrightarrow {} \max_{\widetilde{Y}^*}.  
\end{equation}
Here $\widetilde{Y}^* \subset Y^*$ is a set where the value of the supremum is not equal to~$+\infty.$

We will use the following notations: 
$u$ is the exact solution of problem~\eqref{dka_primal_problem}, $v$~is an~approximate solution of~\eqref{dka_primal_problem}; $p^*$~and~$y^*$ are correspondingly the exact and an approximate solution of problem~\eqref{dka_dual_problem}. If $u$ and $p^*$ exist and are unique, then they are connected via the relation $p^*=\mathcal A \Lambda u.$

\begin{lemma}[Deviation Identity]\label{dka_dev_id}
	For any function $v \in \mathbb K$ and for any $y^* \in \widetilde{Y}^*,$ the~following equality holds:
	\begin{equation}\label{dka_main}
		\begin{aligned}
			\frac{1}{2} \| \Lambda (v-u) \|_{\mathcal A}^2 + 
			\frac{1}{2} \| p^* - y^* \|_{\mathcal A^{-1}}^2 +  \Bigl( \Lambda (v-u), p^*-y^* \Bigr) = \\
			= 	\frac{1}{2} \| \mathcal A  \Lambda v- y^* \|_{\mathcal A^{-1}}^2.
		\end{aligned}
	\end{equation}
\end{lemma}
\begin{remark}
In the right hand side of~\eqref{dka_main} we can also use the norm $\| \Lambda v- \mathcal A^{-1} y^* \|_{\mathcal A}^2.$	
\end{remark}
\begin{remark}
Lemma~\ref{dka_dev_id} is also handy for a ``no-tensor case'', i.e., when $\mathcal A = \mathrm {id}.$ Equivalent norms become standard $L_2$--norms and the deviation identity takes the form
\begin{equation}\label{dka_main_L2}
		\frac{1}{2} \| \Lambda (v-u) \|^2 + 
		\frac{1}{2} \| p^* - y^* \|^2 +  \Bigl( \Lambda (v-u), p^*-y^* \Bigr) 
		= 	\frac{1}{2} \| \Lambda v- y^* \|^2.
\end{equation}
\end{remark}

Now let us discuss the lemma. For linear operators, equality~\eqref{dka_main} is a simple algebraic relation, but it is a good start for further research. In particular, it can be a first step to functional a~posteriori estimates. Such estimates do not need any information on how  approximated solutions were constructed and are both of theoretical and numerical interest, see, e.g.~\cite{dka2020ApRe}. Furthermore, recently they proved to be useful for machine learning PDEs' solvers, see~\cite{dka2024FaOsRu} and references therein. 

The left-hand side of~\eqref{dka_main} shall tend to zero if $v$ is close to $u$ and $y^*$ is close to~$p^*.$ The right-hand side depends only on approximate solutions, thus, is computable. Therefore, we can evaluate the vicinity of approximate solutions to the exact ones. That is more important; we can compare, which approximations give the smallest norm in the r.h.s, so which of them are ``better''. Once more, note that we can make this comparison regardless of methods for approximations' construction.

In Lemma~\ref{dka_dev_id} we used the title ``Deviation Identity'' (DI) instead of the more customary ``Error Identity'' (EI). The last one usually refers to equality in a special form, with so--called compound functionals $\mathcal D_G$ and $\mathcal D_F$ (see, e.g., relation~(2.7) in~\cite{dka2020ApRe}). Although these two identities are very closely connected, and DI can be derived from the EI, we use different names to distinct one from another. In~\cite{dka2021Re} one can find some useful considerations about the deviation identities for problems similar to~\eqref{dka_primal_problem}.

The relation~\eqref{dka_main} itself holds not only for $v \in \mathbb K$ and $y^* \in \widetilde{Y}^*$ but for any $v \in V$ and any $y^* \in Y^*$ (see also discussion in~\cite{dka2023Bes} for the biharmonic case). Moreover, it is valid for any pair $u$ and $p^*$ such that $p^*=\mathcal A \Lambda u.$
So why do we require $\mathbb K$ and $\widetilde{Y}^*$ in the formulation of  Lemma~\ref{dka_dev_id}? For $v,$ the answer is obvious: we are interested only those functions for whom limitations describing the subset $\mathbb K$ hold. We will discuss the reason for having condition 
\begin{equation}\label{dka_infty}
\sup\limits_{v \in \mathbb K} \biggl\{ - \left( y^*, \Lambda v \right) + (f,v) \biggr\} \ne +\infty
\end{equation}
later. For now, let's say that it provides us the ``better situation''.

The most important and the most complicated part of equality~\eqref{dka_main} is the scalar product  $\Bigl( \Lambda (v-u), p^*-y^* \Bigr).$ It contains all the data of our problem: the operator $\Lambda,$ the function $f,$ and the minimization set $\mathbb K$ (and, luckily, it does not contain the tensor $\mathcal A).$ From now on we will denote this scalar product as  $M_{\mathbb K}.$ In the current form, it is not very convenient from the practical point of view. It contains two differences, both for primal and dual variables. Thus, it is inconvenient to estimate the scalar product in the general case. Also,  we don't know a priori the sign of this product or its parts. Nevertheless, for particular cases, the strategy for calculating  $M_{\mathbb K}$ exists. Below, we will show it for the biharmonic obstacle problem. Without loss of generality, we will examine the case where  $\mathcal A = \mathrm {id}.$

\section{Biharmonic obstacle problem}
In this section, we consider $V=\bigl\{ \nu \in H^2 (\Omega) : \nu \bigl|_{\partial \Omega} =\frac{\partial \nu}{\partial n} \bigl|_{\partial \Omega}=0 \bigr\}$ and  $\Lambda = \nabla \nabla.$ The operator $\nabla \nabla$ satisfies the coercivity condition due to the Friedrichs--type inequality. Note that the problem is biharmonic because of the squared norm in~\eqref{dka_primal_problem}. 
Further, $Y=M^{n \times n}_{sym} (L_2( \Omega))$ is the space of symmetrical matrix--valued functions with the scalar product $\left( q, g  \right)_Y = \int\limits_\Omega q \! : \! g \, dx$ where $q \! : \! g = \sum\limits_{i=1}^n \sum\limits_{j=1}^n q_{ij} g_{ij}.$  

Let $\mathbb K = \left\{ \omega \in V \colon \omega \geqslant \varphi\ \text{a.e. in}\ \Omega \right\}$ where the obstacle $\varphi \in C^2(\overline{\Omega}),$ $\varphi \leqslant 0$  on~$\partial\Omega,$ and $\varphi(x) \le \mathrm{c} (\mathrm{dist}(x, \partial \Omega))^2,$ $\mathrm{c}>0$ in $\Omega.$ The last condition guarantees the non--emptiness of the minimization set (see also~\cite{dka2023Bes}). It is not hard to prove that the set $\mathbb K$ is convex and closed.

Under such conditions, problem~\eqref{dka_primal_problem} can be applied to the elasticity theory, where it describes the equilibrium contact of clamped elastic plates over a rigid obstacle. The existence and uniqueness of the solution for this problem have been proved, see~\cite{dka2020ApRe} and bibliography therein, its properties have been studied. 
It is not any less important that the a priori smoothness of the exact solution are known: 
$u \in H^3_{loc}(\Omega) \cap W^{2, \infty}_{loc}(\Omega)$ and $\nabla \nabla u \in W^{2, \infty}_{loc}(\Omega).$ 
\begin{remark}
	Please note that in the current work no general conditions for existence or uniqueness of the exact solutions are provided. Therefore, we suggest to ensure these facts for a particular problem before using~\eqref{dka_main} and deriving anything from it.
\end{remark} 

In the case of an obstacle problem, the domain $\Omega$ is divided into two subdomains $\Omega_f$ and $\Omega_\varphi.$ In $\Omega_f$ equation $\Delta^2 u=f$ holds. In $\Omega_\varphi$ the minimizer coincides with the obstacle $u=\varphi.$ The boundary between the two subdomains is unknown, i.e. we deal with a free--boundary problem (see~\cite{dka2020ApRe} for more details). This distinction for subdomains and conditions will play an essential role when we calculate the scalar product $M_{\mathbb K}.$

Firs of all, we will quite naturally assume that $\mathrm{divDiv} \;y^* \in L_2(\Omega).$	Then we will:
\begin{enumerate}
\item integrate $M_{\mathbb K}=\Bigl( \nabla \nabla (v-u), p^*-y^* \Bigr)$ by parts tacking into account the homogeneous boundary restrictions,
\item recall different conditions in subdomains $\Omega_f$ and $\Omega_\varphi,$
\item take care of the situation on the interface between the subdomains.
\end{enumerate}
These actions will lead to the following result:
\begin{multline*}
	M_{\mathbb K} =-\int\limits_{\partial \Omega_\varphi} (v-u) \left[  \mathrm{Div}p^* \cdot e_{\Omega_\varphi} \right]  dS 
	+ \int\limits_{\Omega_\varphi} (v-\varphi) \left( \mathrm{divDiv}p^* - f \right) dx + \\
	+ \int\limits_{\Omega_f} (u-\varphi) \left(  \mathrm{divDiv}y^* -f \right) dx 
	- \int\limits_\Omega (v- \varphi ) \left( \mathrm{divDiv}y^* -f \right) dx. 
\end{multline*}
Here $[N]=N(\Omega_\varphi)\bigl|_{\partial \Omega_\varphi} - N(\Omega_f)\bigl|_{\partial \Omega_\varphi}$ denotes the jump of $N$ across the free boundary $\partial \Omega_\varphi,$ and $e_{\Omega_\varphi}$ is an outer unit normal to $\partial \Omega_\varphi.$ 

Now in accordance with~\cite{dka2020ApRe}, let us denote:
\begin{equation*}
	\begin{aligned}
	\mu_\varphi(v)&=-\int\limits_{\partial \Omega_\varphi} (v-u) \left[  \mathrm{Div}p^* \cdot e_{\Omega_\varphi} \right]  dS 
	+ \int\limits_{\Omega_\varphi} (v-\varphi) \left( \mathrm{divDiv}p^* - f \right) dx,\\
	\mu^*_\varphi(y^*)&=\int\limits_{\Omega_f} (u-\varphi) \left(  \mathrm{divDiv}y^* -f \right) dx. 
	\end{aligned}
\end{equation*}
Then, for the biharmonic obstacle problem, we get
\begin{multline}\label{dka_main_bih}
	\frac{1}{2} \| \nabla \nabla (v-u) \|^2 + \mu_\varphi(v) +
	\frac{1}{2} \| p^* - y^* \|^2 +  \mu^*_\varphi(y^*) = \\
	= 	\frac{1}{2} \| \nabla \nabla v- y^* \|^2 + \int\limits_\Omega (\varphi -v) \left( f - \mathrm{divDiv}y^* \right) dx. 
\end{multline}
This relation is exactly the a posteriori identity~(2.24) from~\cite{dka2020ApRe}, but obtained with less calculations.

Now let's go back to~\eqref{dka_infty}. In~\cite{dka2020ApRe} it is shown that for the biharmonic obstacle problem, this condition holds if and only if 
$f - \mathrm{divDiv}\;y^* \le 0$ a.e. in $\Omega.$
The remarkable thing is that under such restriction all members in~\eqref{dka_main_bih} are non--negative! Otherwise, the signs of $\mu_\varphi(v),$ $\mu^*_\varphi(y^*)$ and the last integral in the r.h.s of~\eqref{dka_main_bih} may be undefined (see also~\cite{dka2023Bes}). Thereby, for relation~\eqref{dka_main_bih} we assume $$\widetilde{Y}^*=\bigl\{ y^* \in M^{n \times n}_{sym} (L_2( \Omega)) \! : \mathrm{divDiv}\; y^* \in L_2(\Omega) \ \mbox{and} \ f- \mathrm{divDiv}\;y^* \le 0 \bigr\}.$$ 
\begin{remark}  In $\widetilde{Y}^*$ let us add conditions $y^*=0$ and $\mathrm{Div} \;y^*=0$ on $\partial \Omega.$ Then we can consider space $V$ equipped with nonhomogeneous boundary conditions. In this case we will obtain equality~\eqref{dka_main_bih} with precisely the same terms as for the homogeneous problem.   
\end{remark}

\begin{acknowledgement}
This research was sponsored by the Russian Science Foundation, project No.~\mbox{24-11-00073.} 
 The author expresses her gratitude to the Yerevan State University (Yerevan, Armenia) where the results were presented for the first time.
\end{acknowledgement}





\begin{thebibliography} {99}
	
\bibitem{dka2020ApRe}
Apushkinskaya, D.E., Repin, S.I.:
Biharmonic obstacle problem: guaranteed and computable error bounds for approximate solutions.
Computational Mathematics and Mathematical Physics
\textbf{60}(11), 1823--1838 (2020)

\bibitem{dka2024FaOsRu}
Fanaskov, V., Oseledets, I., Rudikov, A.: Neural functional a posteriori error estimates (2024).
DOI: 10.48550/arXiv.2402.05585

\bibitem{dka2021Re}
Repin, S.I.: Identity for deviations from the exact solution of the problem $\Lambda^*\mathcal{A}\Lambda u+l=0$ and its implications.
Computational Mathematics and Mathematical Physics
\textbf{61}(12), 1943--1965 (2021)

\bibitem{dka2023Bes}
Besov, K.O.: Integral identity and estimate of the deviation of approximate solutions of a biharmonic obstacle problem.
Computational Mathematics and Mathematical Physics
\textbf{63}(3),  333--336 (2023)

	
	
	
	
	
	
	

	
	
\end{thebibliography}
\end{document}